\documentstyle{article}
\def\Spec{{ \mbox{Spec} }}
\def\Z{{ {\bf Z} }}
\def\S{{ \Sigma }}
\def\deg{{ \mbox{deg} }} 
\def\Spec{{ \mbox{Spec} }}
\def\1ox{{ \Omega^1_{\scriptstyle{X}} }}
\def\2ox{{ \Omega^2_{\scriptstyle{X}} }}

\def\ok1{{ \Omega^1_K }}
\def\ok2{{ \Omega^2_K }}
\def\Om{{ \Omega }}
\def\lom{{ \Omega (\mbox{log}D) }} 
\def\om{{ \omega  }}
\def\O{{ {\cal O} }} 

\def\ra{{ \rightarrow }}
\def\hra{{ \hookrightarrow }}
\def\da{{ \downarrow }}

\def\M{{ \overline{{\cal M}}_g }}
\def\1M{{ \overline{{\cal M}}_{g,1} }}

\def\hra{{ \hookrightarrow }}

\def\Ra{{ {\Rightarrow} }}
\def\C{{ {\bf C} }}
\def\Q{{ {\bf Q} }}
\def\P{{ {\bf P} }}

\newtheorem{thm}{Theorem}

\title{Height inequalities and canonical class inequalities}
\author{Minhyong Kim}
\date{12/07/97}
\begin{document}
\maketitle
\section{}

This lecture is concerned with apriori bounds on the
size of solutions to Diophantine equations.
We will 
present a rather random collection of results   with the
intention to provide a flavor of the subject, with
only a glimpse of something like a unifying principle.
In the next lecture, we will make a better attempt to
give coherence to our exposition, by describing the
relation to {\em canonical class inequalities} from
the theory of algebraic surface. The beginning of the
lecture is intended to be entirely popular, while the
technical language will grow denser towards the end.

It might be best to begin with an example from ordinary
number theory  which is likely to be more familiar
 than the arithmetic of function-fields, although it
is  with the latter that we
will be most concerned, and for which there are definite results
of a general nature.
Consider the equation
$$x^3+y^3=1729$$
to which we seek solutions in integral coordinates. Ramanujan
observed the existence of the solutions (9,10) and (1,12), and
noted that these were the only ones. How might one go 
about proving this assertion? Here is one method:
The left hand side of the equation factors into
$$(x+y)(x^2-xy+y^2)$$
so one notes that for any solution $(x,y)$, the quantity
$x^2-xy+y^2$ must divide 1729. In particular,
$|x^2-xy+y^2|\leq 1729$. But $x^2-xy+y^2=(x-y/2)^2+3y^2/4$ so
we get $$3y^2/4\leq 1729\Ra |y|\leq 2\sqrt{1729/3}\approx 48$$
and by symmetry, we have the same bound for $|x|$.
Thus, we arrive at a brute force algorithm for finding
all the solutions, namely, by enumerating all the
possibilities up to the given bound, and observing that the
pair above are the only ones. In general, this method shows that
the integer solutions to
$$x^3+y^3=m$$
satisfy the bounds $|x|,|y|\leq 2\sqrt{|m|/3}$. This
bound can be used to verify another of Ramanujan's claim
from the same anecdote, that 1729 is the smallest natural
number which can be written as a sum of cubes in two different
ways.

Finding solutions to Diophantine equations such
as the above, or such as $x^n+y^n=z^n$ is almost always a non-trivial
problem, and simple resolutions are an extreme rarity.
Also, if one is looking for rational solutions in addition
to integral ones, the difficulties multiply, as may be
seen even in the example above: $(20760/1727,-3457/1727)$ is
also a rational solution, and in fact, the equation has infinitely
many rational solutions.
At present one can safely say that the general problem of finding
all rational solutions to a two-variable Diophantine equation
has a satisfactory {\em systematic}
resolution only for equations of genus zero and some
families of genus-one equations,
in addition to scattered examples (e.g., Fermat's equation,
some families of modular curves, certain hyperelliptic equations).
On the other hand, for equations of genus at least two,
Faltings' theorem assures us of the finiteness of solutions,
so it makes sense to ask the question:
Given an (irreducible) equation $f(x,y)=0$ with rational
coefficients and genus at least 2, can one find an
apriori bound for the size of its rational solutions?
Here size of a solution $(x,y)$ may be
defined as follows: Write $(x,y)=(p/r,q/r)$ with $(p,q,r)=1$.
Then the size or {\em height} of $(x,y)$ is defined to be
sup$(|p|,|q|,|r|)$. Thus, giving a bound for the height
ensures that the brute force search is a finite algorithm.

This question is known as the {\em effective Mordell conjecture}.
While it is increasingly clear that this conjecture sits
inside a coherent network of arithmetic-geometric investigations
(Vojta's conjectures, ABC conjectures, Lang's conjectures, etc.),
it is also fair to say that the problem is very
difficult, and there are no positive results
at present. A goal of these two lectures will be to
convince the audience, by bringing out
the underlying structure of the conjecture, that it is nevertheless
reasonable to believe in its validity.

For this, we are motivated to study the same question in the
context of Diophantine equations over function fields,
which have traditionally been a source of inspiration and
techniques that often end up influencing in a decisive
way the corresponding problem over number fields. Furthermore,
the function field set-up has given rise to many 
questions of algebraic geometry interesting in their
own right, independent of their provenance in
number theory. 

So we move on to some examples in that direction.
 Let $f(x,y,t)$ be polynomial in
$\C[x,y,t]$ and suppose we are interested in finding all
rational function solutions to the equation
$$f(x,y,t)=0$$
where by solution, we mean pairs $x(t), y(t) \in \C(t)$
such that
$f(x(t),y(t),t)=0$
as a function of $t$. For example,  the equation
$$y^3=x^4-6tx^3+11t^2x^2-6t^3x$$
 has the solution $(t,0)$ and the equation
$$(t^4+t)y^3=(t^3+1)x^4+tx^3-t^4$$
has the solution $(t,t)$.
The well-known analogy here is between the field $\Q$ of
rational numbers and the field $\C(t)$ of rational functions
over the complex numbers. There have been many results
that give apriori bounds for the degrees of solutions,
of which we start by mentioning one, in simplified form.
Let $f(x,y,t)$ be an irreducible polynomial and
let $X$ be the closure in $\P^2\times \P^1$ of the zero set
$f(x,y,t)=0$ in $\C^3=\C^2\times \C^1$. Thus, $X$is
equipped with a map $\pi$ to $\P^1$ and can be viewed as a projective
surface fibered over $\P^1$, or a projective family of
curves parametrized by $\P^1$, where for each $t\in \P^1$,
$X_t$ is the  closure of $f(x,y,t)=0$in
$\P^2$. We assume that
$X$ is a non-singular surface, and that  the general fiber of
$\pi:X\ra \P^1$, that is, (the projective plane curve corresponding
to) $f(x,y,t)=0$ for a general value of $t$, is connected and smooth.
 Further assume
that the $(x,y)$-degree $d$ of $f$ is at least 4
and let $e$ be the $t$-degree of the polynomial
$f$. We define a two invariants $s, k$ by  looking at
each curve $f(x,y,t)=0$ for fixed $t$ separately, and noting that
only finitely many are singular, while the rest are smooth of genus
$g=(d-1)(d-2)/2$. Denote the number of singular fibers by $s$.
The singular fibers may contain some components
that are rational curves, and these can be found readily. Denote by
$k$ the total number of rational components obtained thus.

All the quantities just defined can be easily
computed if $f$ is given.
We will measure the size of rational solutions as follows:
If $(x(t),y(t))$ is a pair of rational functions, then
write $x=p(t)/r(t),y=q(t)/r(t)$ where $(p,q,r)=1$, and define
$h(x,y)=\mbox{supdeg}(p,q,r)$. Among many results
on apriori bounds, we mention one, which is in many senses
the sharpest \cite{T}:

\begin{thm}(S.-L. Tan)
With the assumptions above, one has the following bound
for the height of solutions in $\C(t)$ to the 
equation $f(x,y,t)=0$:

$$h(x,y)\leq\frac{(d^2-3d+1)(s-1)+k}{d-3}$$

\end{thm}

This allows us in  principle to find all solutions
to the equation by putting $N$ to be the right-hand side
of the inequality and doing a subtitution
$$(x,y)=(\frac{p_0+p_1t+\cdots +p_Nt^N}{r_0+r_1t+\cdots+ r_Nt^N},
\frac{q_o+q_1t+\cdots+q_Nt^N}{r_0+r_1t+\cdots+ r_Nt^N},$$
to find all solutions via Groebner basis techniques.

Where does this kind of inequality come from and why does one
get so much stronger results than the number field
case? One answer might be that
 the algebraic geometry of surfaces gives us powerful
tools for dealing with  problems of this kind. More
specifically,  the
main technical advantage over the arithmetic version of the
problem is the use of differential forms on curves and surfaces. 
To explain this, let us recall the set-up for translating the Diophantine
problem into algebraic geometry.

Given a polynomial $f$ as above, we can view it as
defining a curve $X_{\eta}$ of genus at least two
over $\eta= \Spec(\C(t))$ which, in turn,
is the generic point of $B=\P^1$. Thus, we can thicken
the curve into a family
$$\begin{array}{ccc}
X_{\eta } &\hra & X \\
\downarrow & & \downarrow \scriptstyle{f}\\
\eta &\hra & B
\end{array}$$
that is, an algebraic surface fibered over $B$. 
There is actually a canonical way of constructing
this fibration, which is called the regular, relatively
minimal model.
Solutions of
the original equations can then be interpreted as diagrams
$$\begin{array}{ccc}
& &X \\
 &\stackrel{P}{ \nearrow} &\da \\
 B&\stackrel{Id}{\ra} & B
\end{array} $$
or {\em sections} of the fibration, in the terminology
of algebraic geometry. We will also be more generally interested
in the multi-sections of the fibration, that is, diagrams
$$\begin{array}{ccc}
& &X \\
 &\stackrel{P}{ \nearrow} &\da \\
 T&\ra & B
\end{array} $$
where $T\ra B$ is a finite covering,
corresponding to solutions of the original equation with
coefficients in some finite extension of the function field
of $B$, i.e., an algebraic solution.
Then given a projective embedding $X\hra \P^n$,
one can view $P$ as being a map $P:T\ra \P^n$ and use the geometry
of $X$ to give bounds on the degree of this map. Typically,
the results are phrased in terms of bounds for 
$$h(P):=\frac{1}{[T:B]}\mbox{deg}(P^*\omega_{X/B}),$$
where $\om_X:=K_X\otimes f^*K_B^{-1}$ is the relative dualizing
sheaf. Being intrinsic to the surface, it is easier to work with,
and it is near enough to being ample to actually give height bounds.

For example, in the case of a surface in $\P^2\times \P^1$
 of
bidegree $(d,e)$ fibered over $\P^1$,
we get
$\om=\O (d-3,e)$
Also, when considering algebraic points,
an important invariant  is the relative discriminant
$d(P):=(2g_T-2)/[T:B]$, which measures the ramification in the
map $T\ra B$.

Then the height inequalities mentioned above follows from a
geometric result,

\begin{thm}(S.-L. Tan)
Let $X$ be a projective smooth surface, $B$ a projective smooth
curve. Let $f:X\ra B$ such that the generic fiber is
smooth and connected
 of genus at least two and the fibration is relatively
minimal. Then for every section $P$, we have
$$h(P)\leq (2g-1)(d(P)+3s)-\om^2.$$
where $s$ refers to the number of singular fibers in the
fibration $X/B$.
\end{thm}
This result is an immediate consequence of a logarithmic
version of the Miyaoka-Yau inequality
for algebraic surfaces of general type, due
to Sakai and Miyaoka \cite{M}, \cite{S}. As such, it uses
the differentials $\Om_X$ (in fact, even the log differentials)
on the surface $X$. Once one has degree bounds of this
sort, one knows that the possible rational points run through
a finite algebraic family of divisor. On the other hand,
if the fibration is  not `isotrivial', that is, is not
birational to a product of curves even after a finite
base extension, any rational point $P$ in the above
sense satisfies an intersection-theoretic
inequality $<P.P><0$. Thus, each $P$ is rigid, so that the
boundedness of the family actually implies finiteness. So
inequalities of the above form are always much stronger
statements than general finiteness results, that is,
Manin's theorem.

A class of results that uses the differentiation on the
base $B$, which is often cited as being the distinguishing
characteristic of function fields, 
involves the {\em Kodaira-Spencer map} of the fibration
which measures the  variation of the family.

There are many different kinds of K-S maps that one can associate
with a family of curves, but the one that seems to have the
most relation to the Diophantine properties seems to be
the K-S map measuring the variation of the Jacobians (that is, the
period integrals).
We recall the definitions. If we have a fibration $X/B$ as above,
choose an open set $U\subset B$ such that $f:X_U\ra U$ is smooth.
We get an exact sequence
$$0\ra f^*\Om_U\ra \Om_X \ra \Om_{X/U}\ra 0$$
from which we get  a map
$$f_*(\Om_{X/U})\ra R^1f_*(\O_X)\otimes \Om_U$$
as the boundary map for the higher direct images.
Passing to the limit over open sets $U$, we get
$$\mbox{K-S}:H^0(\Om_{X_{\eta}})\ra H^1(\O_{X_{\eta}})\otimes \Om_{\eta}$$

This map is actually the cotangent map of the map to the moduli
stack of principally polarized
abelian varieties given by Jacobian of the curve $X_\eta/\eta$.
It is a fact that the generic curve over a function field
$\eta$ will have a K-S map of full rank, that is, will be
maximally varying.
The Kodaira-Spencer map already occurs in the original
proof that Manin \cite{Ma} gave of the geometric Mordell conjecture.
This is because the K-S map is a `component' of the
Gauss-Manin connection. For example, the assumption that
K-S map is of full rank is equivalent to the assumption
that the Gauss-Manin connection is generated by the global
relative differentials.

We state a few more height inequalities in this vein.
\begin{thm}(Moriwaki \cite{Mo} complemented by  \cite{K2})
Let $X$ and $B$ be as above and suppose the K-S map of
$X/B$ has full rank. Then for any algebraic point $P$, we
have
$$h(P)\leq 4d(P)+4c_2(X)-c_1(X)^2-4(g_B-1).$$
\end{thm}

This inequality is striking in that the dependence on
the discriminant $d(P)$ is very weak.
Serge Lang \cite{L} has posed the question of finding the
best possible constants for which one can prove
inequalities of the form
$$h(P)\leq Ad(P)+B$$
for algebraic points $P$,
where $A$ and $B$ should depend on the surface $X$.
The results above are all motivated by this question, which
is related to the algebraic geometric investigation of the
geography of surfaces. In the theorem above, we see that
one can take $A$ independent of the surface $X$ for a generic
family.

Another result that should also be mentioned in this regard
is the \cite{V}
\begin{thm}(Vojta)
With $X$ and $B$ as above, we have
$$h(P)\leq (2+\epsilon)d(P)+O(1)$$
where O(1) depends on $X$ and $\epsilon>0$, but not $P$.
\end{thm}

This results is a function-field version of Vojta's conjecture on
algebraic points. This result can be extended to characteristic
$p$ under the assumption again that the K-S map is
of full rank (\cite{K1}), which is a stronger assumption in
positive characteristic than in characteristic zero (because
of the possibility of inseparability). In relation to
Lang's question, one has made the $A$ term (dependence on $d(P)$)
very small at the expense of making $B$ inexplicit.

Diophantine equations over function-fields of
 characteristic p are different again from the case of complex
function fields. One advantage over complex function-fields,
as
far as the closeness to number fields goes, is the fact that
one can have finite residue-fields. On the other hand, because
of the presence of inseparable maps, many difficulties have a flavor
entirely unique to this domain. Here is one general result:
\begin{thm}(\cite{K1})
Suppose $X/B$ is not isotrivial. Then
$$h(P)\leq p^e(2g-2)d(p)+O(\sqrt{h(p)})$$
where $p^e$ is the degree of inseparability of
the classifying map, that is the map $B\ra {\cal M}_g$,
for the family $X/B$.
\end{thm}
The proof of this result also goes through the K-S map.
Voloch has shown how  a new proof of the Riemann
hypothesis for curves over finite fields follows from this
inequality. That is, the square root of the height inside
the big-Oh is the same one that occurs in estimates for
numbers of points on curves. So there are  interesting
arithmetical connections which have not at all been investigated in
detail.

Because of the Frobenius maps, the assumption of non-isotriviality
is important in characteristic p.
That is, for example, if $X=B\times B$, and
$B$ is defined over a finite field, then some power of
the Frobenius map will give a morphism $F^n:B\ra B$, so
morphisms $F^{kn}:B\ra B$ for each $k>0$. Thus we get infinitely
many rational points of $X=B\times B \ra B$. There are
various more complicated `twisted' situations where
similar phenomena can arise. Purely inseparable
maps are pervasive objects in characteristic
p geometry, which must always be dealt which separately.
For example, consider the following amusing problem: Given a fibration
$X\ra B$, how many diagrams of the following sort can one have?
$$\begin{array}{ccc}
& &X \\
 &\stackrel{P}{ \nearrow} &\da \\
 B&\stackrel{F^n}{\ra} & B
\end{array} $$
An elementary way of stating the question is to consider a
polynomial $f(X,Y)$ with coefficients in a characteristic $p$
function field $K$, and to start `twisting' the equation with
the Frobenius map of $K$, that, raising the coefficients of
$f$ to the $p$-th power. Denote by $f^{(n)}$ the polynomial
twisted $n$ times. Notice that  if $(X,Y)$ is a solution of
$f(X,Y)=0$, then $(X^p,Y^p)$ is a solution of $f^{(1)}=0$.
That is, old solutions can be twisted to give solutions
of the twisted equation. The question being posed, then,
is whether or not we can keep getting {\em new} solutions
as we twist $f$. These can be view as solutions to $f$
in the purely inseparable tower $K^{1/p^{\infty}}$.
There is a sense in which `purely-inseparable points' of
this sort should not be that different from rational points since
the map $F^n:B\ra B$ is just a homeomorphism.
This sentiment is captured by the following result,
which is more natural to 
express directly in terms of a curve over a function field
rather than a surface fibered over a curve:
\begin{thm}( \cite{K3})
Suppose $K$ is a function field of chracteristic $p>0$,
and suppose $C/K$ is a smooth curve of genus at least 2.
Finally, suppose $C$ is not isotrivial, that is, $C$ does not
become isomorphic to a curve defined over a finite field after
some finite base-extension. Then $C(K^{1/p^{\infty}})$ is
finite, where $K^{1/p^{\infty}}$ is the tower of all purely
inseparable field extensions of $K$.
In fact,  if we assume that the curve
has semi-stable reduction over the smooth model $B$ of $K$, then
the purely-inseparable points satisfy a height inequality
$$h(P)\leq (2g_B-2+s)$$
where $s$ is the number of singular fibers in the semi-stable reduction.

\end{thm}

All the results mentioned use  differential forms
in one or more serious ways, posing for arithmeticians
the problem of finding a number field substitute for
differentiation. This problem has  occupied the efforts
many number-theorists for the last few decades.

\section{}

 The point of view we will emphasize here for understanding
the height inequalities mentioned in the previous lecture
is that they are all instances of {\em canonical class inequalities}.
That is, if one didn't know of it beforehand, one would expect
such height inequalities to exist, simply because of the pervasiveness
of canonical class inequalities in algebraic geometry.

To elaborate a bit, a 
class of geometric objects will typically come with a set
of numerical invariants which are often related by various
natural formulas. As an example, one may consider Noether's
formula for algebraic surfaces. In contrast to such equalities,
one also has general {\em inequalities} between numerical invariants
which are of aid in classification theory, and which are often 
more elusive than the equalities.

For example, for a surface of general type $X$ over a field of characteristic
zero,  we have \cite{M}, 
$$c_1(X)^2\leq 3c_2(X),$$
known as the Miyaoka-Yau inequality,
which arose  from a long sequence of results due to Van-de-Ven, Bogomolov,
Miyaoka and Yau \cite{BPV}.

That is, one attempts to classify surfaces of general type according to 
the pair of numerical invariants $(c_1^2,c_2)$. Those
surfaces with a given pair of invariants form a moduli
space in a natural way (Gieseker). However, it is still
an ongoing project in classification theory to determine
which of these schemes are non-empty, that is, which pairs
of integers $(n,m)$ actually occur as $(c_1^2,c_2)$. This
problem is often referred to as one of determining the
`geography' of the surfaces of general type.
If one  plots the possible values in the plane, the M-Y inequality
says that all the values lie above the line $c_2=3c_1^2$.
For reference, we quote from the book of Barth-Peters-Van-de-Ven \cite{BPV}
 some of the other facts  known
about the distribution of the $(c_1^2,c_2)$'s for surfaces of general type,
i.e. some other canonical class inequalities:

-Clearly, $c_1^2+c_2\cong 0$ (mod 12);

-$c_1^2>0$ and $c_2>0$;

-$5c_1^2-c_2+36\geq 0$ if $c_1^2$ is even;

 $5c_1^2-c_2+30\geq 0$ if $c_1^2$ is odd; (Noether's inequalities).

The point we wish to emphasize is that such equalities and
inequalities between numerical invariants are pervasive in
geometry.

Another class of canonical class inequalities is associated to
the geometry of the moduli space $\M$ of stable curves of
genus $g\geq 2$. On $\M$, there are several
 natural divisors classes. One is
 the class $\Delta$ of the boundary, corresponding to
the singular curves, and two more classes arise from
the universal curve $f:{\cal C}\ra \M$.  One is the
first chern class of the Hodge bundle, $\lambda=c_1(\pi_*(\om))$,
where $\om$ is the relative dualizing sheaf of the universal curve,
and the other is the direct image $k=\pi_*(\om.\om)$ of the
self-intersection of $\om$ which gives
a divisor on $\M$. It is well known that these divisors (and the
universal curve) are only defined on the corresponding
moduli stack, but there is also a well-known isomorphism of
rational chow groups, which is where our classes lie. They
are also ${\bf Q}-$Cartier divisors, and hence can be intersected
with curves.
From these classes one can construct numerical invariants for
any family of stable curves parmetrized by a projective
curve. That is,
if $f:X\ra B$ is such a family, then we get a map
$B\ra \M$, and so numbers $\delta_X, \lambda_X$ and $k_X$,
obtained by pulling back the ${\bf Q}-$Cartier divisors
above and taking degree. 
These can be interpreted directly in terms of the family itself:
$k_X=<\om_{X/B}^2>$ and $\lambda_X=\mbox{deg}f_*\om_{X/B}$.
$\delta_X$ is perhaps best explained by taking a minimal
desingularization $p:X'\ra X$ and taking the number of
singularities in the fibers of $X' \ra B$. Note that here,
we are taking into account the fact that the original $X$
may have been singular as a surface, not just have
singularities in the fibers. $X'$ is usually refered to as
the regular semi-stable model of $X$. In any case, this
way we arrive at natural invariants of a family $X/B$.

One general inequality due to Cornalba-Harris \cite{CH} and Xiao {X} is
that for any stable family, we have:
$$(1-1/g)\delta \leq (2+1/g) \om^2$$
That is, $\om^2$ is bounded below in term of the
number of singularities. One the  other hand, the
invariants for $X'$ as a surface can be related to the
invariants of the family $X/B$, where we assume that
$B$ also has genus at least two, so that $X'$ will be
of general type. The relation arises from the
exact sequence
$$0\ra f^*\Om_B \ra \Om_X' \ra I_Z\om_{X'} \ra 0$$
where $I_Z$ is the sheaf of ideals for the singular
points of the fiber in $X'/B$. One also has
$p^*\om_X=\om_{X'}$ ( $X'\ra X$ is a crepant resolution),
so that $\om_X^2=\om_{X'}^2$. 
If we put these facts together,
we find that $c_1^2(X')=(\om_X+K_B)^2$ while $c_2(X')=\om_X.K_B +\delta$,
where we will abuse notation a bit by basically confusing $X$ and $X'$.
So Noether's formula, for example, reads
$$12\deg f_*(\om)=\om^2+\delta,$$ 
while Noether's inequality translates into:
$$\delta\leq 5\om^2+9(2g-2)(2g_B-2)+36$$
which is similar to the C-H-X inequality, but weaker.
 This version of Noether's formula is
also responsible for the fact that our three
invariants are not independent, so that any other relation involving
 two of them, such as the preceding inequality, actually has obvious
consequences for the remaining one.
 How does the M-Y inequality apply to families?
Simple arithmetic yields
$$\om^2\leq (2g-2)(2g_B-2)+3\delta ,$$
an {\em upper} bound for $\omega^2$ in terms of delta.
To summarize, canonical class inequalities in two different
settings, that of surfaces of general type, and that of
families of curves, are related in a natural way.  We just
described the relation in one direction, that is, from surfaces
to curves, but it is not hard to work out implications that
the C-H-X inequality, for example, has for Chern numbers of surfaces, by
choosing a canonical model and a Lefshetz fibration.

One last result in this vein we mention is an improvement on
the C-H-X inequality due to Eisenbud, Harris, and Mumford,
which holds for generic families of stable curves:
$$\delta \leq (1+o(1/g))\om^2.$$

What do these results have to do with height inequalities?
And why was it stated above that height inequalities
are natural things to
expect?
The point is that height inequalities for rational points
are nothing but canonical class inequalities for {\em families of
pointed curves}. That is, they relate the invariants
arising from $\1M$, the moduli space of stable pointed curves.
For a pointed stable curve

$$\begin{array}{ccc}
& &X \\
 &\stackrel{P}{ \nearrow} &\da \scriptstyle{f} \\
 B&\stackrel{Id}{\ra} & B
\end{array} $$
the natural analogue of the sheaf $\om$ is the sheaf
$\om (P)$. One can give many justification for this
analogy, but the most important one is perhaps that
$f_*(\om^{\otimes2} (P))$ at the generic point
is the cotangent space to
$\1M$ at the point corresponding to the
pointed curve: that is, $\om(P)$ governs the variation of
a pointed curve in the same way that $\om$ governs the variation
of a curve.
Thus, one gets the new invariants, deg$f_*(\om(P))$,
$(\om(P).\om(P))$, and others which need not concern us
at the moment. Now
$$(\om(P).\om(P))= \om^2+2(\om.P)+P^2=\om^2+(\om.P)$$
since $P^2=-(\om.P)$ by the adjunction formula.
So the height of the point $P$ is the most natural
new invariant that arises  out of the situation when considering pointed
curves rather than curves. 

On the non-fibered side of the story, the relevant theory is
that of {\em log surfaces},
the category of pairs $(X,D)$, consisting
of a surface $X$ and a divisor $D$ on $X$. A point of view which 
has become very current in higher-dimensional classification
theory, and espoused by
Kollar in his Santa Cruz lectures, is that this category of pairs is a 
more natural and fundamental
framework for classification than that of varieties by
themselves, especially in the study of singularities.

Miyaoka and Sakai \cite{S} have also proved a version of the M-Y 
inequality for
log-surfaces. The generalization that is relevant for our
purposes is
 in the situation where
 $(X,D)$ is a pair such that $X$ is smooth and
$D$ has
normal crossings, and concerns the log differentials
 $\lom$. Suppose $c_1(\lom)=K_X+D$ is nef. Then
$$c_1(\lom)^2\leq 3c_2(\lom).$$ 

We outline how this leads to S.-L. Tan's inequality (The following
proof was shown to me by N. Shepherd-Barron):
Suppose $f:X\ra B$ is as in the previous lecture,
and we assume for convenience that the family is
semi-stable. Let $P\subset X$ be the image of
a section and let $\S$ be the sum of the singular fibers,
and $S=f(\S)\subset B$ the points lying below the singularities.
Then $D=P+\S$ is a normal crossing divisor, so we
consider the log differentials $\lom$. There is an
exact sequence
$$0\ra f^*\Om_B(S)\ra \lom \ra \om(P) \ra 0$$
arising from the fact that the relative dualizing sheaf
$\om$ can be identified with the relative log differentials.
Thus $c_1(\lom)=\om+P+K_B+\S$ so that 
$$\begin{array}{rcl}
c_1^2(\lom)&=&\om^2+2<\om.P>+<P^2>+2<\om+P.K_B+\S>\\
 &=&\om^2+<\om.p>+2<\om+P.K_B+\S>
\end{array}$$
while 
$$c_2(\lom)= <\om+P.K_B+\S>.$$So
Miyaoka's inequality reads
$$\om^2+<\om.P>\leq <\om+P.K_B+\S>=(2g-1)(2g_B-2+s)$$,
since the intersection of a horizontal divisor and a divisor
coming form the base is the product of the degrees.
This is exactly the first height inequality we stated.
One can use the lower bound for $\om^2$ given by
C-H-X to eliminate the $\om^2$ in favor of the more
explicit invariant $\delta$.
In any case, this derivation makes clear the point that
a height inequality is nothing but a canonical class inequality.
That is, and this is the main point we wish to convey, 
far from being a result 
of purely
arithmetic interest, height inequalities sit squarely inside
a well-established body of theory from the mainstream of
algebraic geometry.

In the last fifteen years, much effort has been expended
in developing an intersection theory for arithmetic
schemes, that is, projective schemes that are flat over
$\Spec (\Z)$. The inspiration of course has been 
intersection theory for varieties fibered over a curve,
and this has motivated much of the reformulation
of results from surface theory into a fibered setting,
some of which we saw here.
Once we see the results in fibered terms,
they can then serve as a springboard for the
study of arithmetic surface, that is,
families of curves over $\Spec(\Z)$, for the
formulation of conjectures and as well as
for ideas and techniques that can be emulated. 
The arithmetic intersection
theory started by Arakelov, was developed by Faltings \cite{F},
Deligne \cite{D} , Gillet and Soul\'{e} \cite{GS}, 
and many others and an increasing
number of
 results from ordinary intersection theory
are being carried over into this arithmetic setting.
It has come to seem reasonable to believe that
every result for algebraic surfaces has a compelling
arithmetic analogue, whose proof will involve
translations and refinements of the geometric theory.
Faltings proof of Mordell's conjecture was to a large extent
inspired by this viewpoint, and one spectacular application
of the intersection theory was the second proof given by
Vojta, which followed a line of reasoning entirely
analogous to a function field theorem he had proved
earlier. This of course was extended by Faltings to include
finiteness theorems for certain sub-varieties of abelian 
varieties. In all, the benefit of having a clear
understanding of the arithmetic geometry of function fields
has become generally accepted, as well as the belief
in a strong analogy and interplay 
between  geometric and arithmetic theories
of various sorts. Then, having demonstrated the
occurence and application of canonical class inequalities in a
 setting  common  to geometry and
arithmetic, it seems reasonable to hope that
canonical invariants  for arithmetic surfaces will 
satisfy relations of a similar nature. This has already
been demonstrated for a good class of equalities, including
the adjunction formula, Riemann-Roch theorems, and
Noether's formula \cite{F}. Most recently, an important inequality,
the 
positivity of $\om^2$ has been established by E. Ullmo.
 The other interesting inequalities have remained elusive,
in particular,  the one which would give us an effective
Mordell conjecture. However, the remarkably 
coherent structures emphasized here, tying
 together the many different
theories through precise analogies, can hardly allow us to
believe otherwise than that they also remain merely to
be discovered.

\end{document}